\documentclass[14pt]{extarticle}
\usepackage{amsmath,amstext, verbatim,remreset}
\usepackage[T2A]{fontenc}
\usepackage[cp1251]{inputenc}
\usepackage[english]{babel}
\usepackage{amsfonts,amssymb}

\newcommand{\RR}{\mathbb R}
\newcommand{\NN}{\mathbb N}
\newcommand{\ZZ}{\mathbb Z}
\newcommand{\A}{A_\varepsilon}
\usepackage{mathrsfs}
\newtheorem{theorem}{Theorem}

\begin{document}
\begin{Large}
\begin{center}
{\bf V. F. Babenko, O. V. Kovalenko}

{\bf On interpolation and extremal properties of periodic perfect splines}
\end{center}
\end{Large}
\begin{abstract}
{Existing and extremal property of periodic perfect spline, which interpolates given function in the mean were proved.}
\end{abstract}

Denote by $C^m$ ($m \in \ZZ_+$) the space of $m$ times continuously differentiable (continuous when $m=0$) $2\pi$-periodic functions; denote by $L_\infty$ the space of all measurable $2\pi$-periodic functions with finite norm $\|f\| = \|f\|_\infty.$
For $r\in\NN$ denote by $L^r_\infty$ the space of functions $f \in C$ such that $f^{(r-1)}$ is absolutely continuous and $f^{(r)}\in L_\infty$.

For $r\in\NN$ the function $s(t)\in C^{r-1}$ will be called a {\it spline} of order $r$ with knots in the points 
\begin{equation}\label{knots}
t_0<t_1<\ldots<t_n<t_{n+1}:=t_0 + 2\pi,
\end{equation} 
if on each interval $[t_i,t_{i+1}]$, $i=0,1,\dots,n,$ $s(t)$ coincides with restriction to this interval of an algebraic polynomial of order less than or equal to $r$. 

Spline $s(t)$ of the order $r$ with knots in points~\eqref{knots} will be called {\it perfect} if $s^{(r)}(t)=(-1)^i\epsilon$, $\epsilon\in\{1,-1\}$, $t\in (t_{i},t_{i+1})$, $i=0,\dots,n$.

For $r,n\in\NN$ denote by $\Gamma^r_n$ the set of $2\pi$-periodic perfect splines of the order $r$  with less than or equal to $n$ knots.

Perfect splines play important role in exact solution of many extremal problems of approximation theory (see for example~\cite{korn2} and~\cite{korn3}).

The main result of this work is the following theorem.

\begin{theorem}\label{th1}
Let $r,m\in\NN$, $\varphi_k(x):\RR\to\RR$ --- integrable even functions with compact support $[-\varepsilon_k,\varepsilon_k]$, positive on $(-\varepsilon_k,\varepsilon_k)$ and such, that
$$
\int\limits_{-\varepsilon_k}^{\varepsilon_k}{\varphi_k(x)dx}=1,\,k=1,\dots,{2m+1}. 
$$

Let also the numbers $0<x_1<x_2<\ldots<x_{2m+1}<2\pi$ such, that supports $[x_k-\varepsilon_k,x_k+\varepsilon_k]$ of the functions $\varphi_k(x-x_k)$ be pairwise disjoint and are contained in $[0,2\pi)$, $k=1,\dots,{2m+1}$, be given. Then for every function $f\in L_{\infty}^r$ there exist a number $\xi $ and a spline $s(t)\in \xi\Gamma^r_{2m}$ such, that $$\int\limits_{x_k-\varepsilon_k}^{x_k+\varepsilon_k}{\varphi_k(x-x_k)} s(x)dx=\int\limits_{x_k-\varepsilon_k}^{x_k+\varepsilon_k}{\varphi_k(x-x_k)}f(x)dx ,\,{k=1,\dots,2m+1}.$$ Moreover
\begin{equation}\label{extremalProperty}
\|s^{(r)}\|_{\infty}=|\xi|\leq\|f^{(r)}\|_{\infty}.
\end{equation}

\end{theorem}

From theorem~\ref{th1}, in particular, the possibility of simple and multiple interpolation by periodic perfect splines follows.

\begin{theorem}\label{th3} Let $r,m\in\NN$, $0\leq x_1<x_2<\ldots<x_{l}<2\pi$, $k_1,\dots,k_l \leq r-1$ --- integer non-negative numbers such that $\sum_{j=1}^l{(k_j+1)}=2m+1$. Then for every function $f\in L_{\infty}^r$ there exist a number $\xi$ and a spline $s(t)\in \xi\Gamma^r_{2m}$ such that $$s^{(j)}(x_i)=f^{(j)}(x_i),\, j=0,\dots,k_i;\,i=1,\dots,l.$$  Moreover, the extremal property~\eqref{extremalProperty} holds.

In particular, if numbers $0\leq x_1<x_2<\ldots<x_{2m+1}<2\pi$ and a function $f\in L^r_{\infty}$ are given, then there exist a number $\xi$  and a spline $s(t)\in \xi\Gamma^r_{2m}$ such that $$ s(x_i)=f(x_i),\, i=1,\dots,2m+1.$$ Moreover, the extremal property~\eqref{extremalProperty} holds.
\end{theorem}

Theorem~\ref{th3} gives periodic analogue of the results, that were achieved by S.~Karlin~\cite{Karlin1, Karlin2} in non-periodic case. The methods of the proof, that we use, are much simplier than those, that were used by S.~Karlin, and can be used in solution of other problems.

We will show the main steps in the proof of the theorem~\ref{th1}. Theorem~\ref{th3} can be proved with the help of theorem~\ref{th1} readily.

Below we will denote by $*$ the operation of convolution for $2\pi$-periodic functions. Let $f\in C$. Denote by  $\nu(f)$ the number of sign changes of the function $f$ on the period. For $\varepsilon>0$ and $x\in\RR$ denote by 
$$\A(x):=\frac1{2\pi}\sum_{j=-\infty}^{\infty}{\frac{{\rm e}^{ijx}}{ch\varepsilon j}}.$$
Note several properties of the kernel $A_\varepsilon$. If $f\in C$, then $(\A*f)$ is analytical function on real line (see.~\cite{Akhiezer37},\cite{Akhiezer38}; see also~\cite{Nikolskiy46}, \S3); $(\A*f)$ uniformly converges to the function $f$ when $\varepsilon\to 0$ on real line; $\nu(\A*f)\le \nu(f)$ (see for example~\cite{Mairhuber});
$\A$ is even; $\int_\RR{\A (x)dx}=1$. Below we will denote by $B_{r}(t)$, $r\in\NN$, the Bernoulli kernel of order ${r}$ (see for example~\cite{korn2}, \S 3.1).

Set $C_k:=\int\limits_{x_k-\varepsilon_k}^{x_k+\varepsilon_k}{\varphi_k(x-x_k)f(x)}dx$, $k=1,2,\dots,2m+1$. For $\varepsilon > 0$ consider the extremal problem
\begin{equation}\label{extremalProblem}
\int\limits_0^{2\pi}{\left|c+\sum_{k=1}^{2m+1}{c_k(\A*B_r*\varphi_k)(x-x_k)}\right|dx}\to \min 
\end{equation} under conditions
\begin{equation}\label{condition1}
\sum_{k=1}^{2m+1}{c_kC_k}=1,
\end{equation}
\begin{equation}\label{condition2}
\sum_{k=1}^{2m+1}{c_k}=0.
\end{equation}
The convolution with the kernel $A_\varepsilon$ is needed to make the integrated functions analytic. This ensures that the function being minimized is continuously differntiable with respect to parameters $c$ and $c_k$, $k=1,\dots,2m+1$. It is easy to prove that the minimum in the extremal problem~\eqref{extremalProblem} -- \eqref{condition2} exists. 

Consider the Lagrange function:
$$L:=\eta \int\limits_0^{2\pi}{\left|c+\sum_{k=1}^{2m+1}{c_k(\A*B_r*\varphi_k)(x-x_k)}\right|dx}+\lambda\left(\sum_{k=1}^{2m+1}{c_kC_k}-1\right)+\mu \sum_{k=1}^{2m+1}{c_k}.$$

According to method of Lagrange multipliers (see for example~\cite{Metod_Lagr}, \S~2), there exist numbers ${\eta_\varepsilon,\lambda_\varepsilon,\mu_\varepsilon\in\RR}$ such that $\eta_\varepsilon^2+\lambda_\varepsilon^2 + \mu_\varepsilon^2 \neq 0$ (we can count, for example, that $\eta_\varepsilon^2+\lambda_\varepsilon^2 + \mu_\varepsilon^2 =1)$, and $c^\varepsilon,c_1^\varepsilon,c_2^\varepsilon,\dots,c_{2m+1}^\varepsilon\in\RR$, which  fulfills the equalities~\eqref{condition1} and~\eqref{condition2} and such that
\begin{equation*}
\eta_\varepsilon\int\limits_0^{2\pi}{{\rm sgn}\left[c^\varepsilon+\sum_{k=1}^{2m+1}{c_k^\varepsilon(\A*B_r*\varphi_k)(x-x_k)}\right]dx}=0
\end{equation*}
and for $j=1,2,\dots,2m+1$
\begin{equation*}
\eta_\varepsilon\int\limits_0^{2\pi}{{\rm sgn}\left[c^\varepsilon+\sum_{k=1}^{2m+1}{c_k^\varepsilon(\A*B_r*\varphi_k)(x-x_k)}\right](\A*B_r*\varphi_j)(x-x_j) dx}+\lambda_\varepsilon C_j+\mu_\varepsilon=0.
\end{equation*}

Passing to the limit when $\varepsilon\to 0$ we get that there exist numbers $\eta, \lambda, \mu$, $ \eta^2 + \lambda^2 + \mu^2 = 1$, and a function $g(x)$, which is non-zero almost everywhere and such that $\nu(g)\leq 2m$, for which $\eta\int\limits_0^{2\pi}{{\rm sgn}g(x)}dx=0$,
and for $k=1,2,\dots,2m+1$
\begin{equation}\label{5''}
\eta\int\limits_0^{2\pi}{{\rm sgn}g(x)(B_r*\varphi_k)(x-x_k) dx}+\lambda C_k+\mu=0.
\end{equation}
Assumption that $\eta = 0$ or $\lambda = 0$ leads to contradiction. This means that for the function $h(x):=({\rm sgn}g*B_r)(x)$ we have $h^{(r)}(x) ={\rm sgn} g(x)$, and hence $h(x)\in \Gamma^r_{2m}$.

From~\eqref{5''} we get that for $k=1,\dots,2m+1$
$$0=\eta\int\limits_0^{2\pi}{{\rm sgn}g(x)(B_r*\varphi_k)(x-x_k) dx}+\lambda C_k+\mu=$$ $$=(-1)^r\eta ({\rm sgn}g*B_r*\varphi_k)(x_k)+\lambda C_k+\mu=(-1)^r\eta(\varphi_k *h)(x_k)+\lambda C_k+\mu=$$ $$=(-1)^r\eta\int\limits_{x_k-\varepsilon_k}^{x_k+\varepsilon_k}{\varphi_k(x-x_k)}h(x)dx+\lambda C_k+\mu $$(the last equality holds because the functions $\varphi_k(x)$ are even).

This means that the perfect spline $s(x):=-\frac{(-1)^r\eta h(x)+\mu}{\lambda}$ is desired.

We will show that $\|s^{(r)}\|_{\infty}\leq\|f^{(r)}\|_{\infty}$ now. Suppose the contrary. Let $\|s^{(r)}\|_{\infty}>\|f^{(r)}\|_{\infty}$. Set $\delta(x):=s(x)-f(x)$. Then the function$\delta^{(r)}(x)=s^{(r)}(x)-f^{(r)}(x)$ can change the sign only in the knots of  the spline $s(x)$, and hence $\nu(\delta^{(r)})\le 2m$. Note that ${\int_{x_k-\varepsilon_k}^{x_k+\varepsilon_k}{\varphi_k(x-x_k)}\delta(x)dx=0}$, function $\delta(x)$ is non-zero on each interval in $[0;2\pi)$ and $\varphi_k(x)$ are non-negative functions, $k=1,\dots,2m+1$, and hence on each interval $(x_k-\varepsilon_k,x_k+\varepsilon_k)$ the function $\delta(x)$ changes sign. The supports of the functions $\varphi_k(x-x_k)$, $k=1,\dots,2m+1$ are pairwise disjoint, and hence $\nu(\delta)\geq 2m+1$, which lead to contradiction. Theorem is proved.

\end{document}


\begin{Large}
{\bf УДК 517.5}
\begin{center}
{\bf В. Ф. Бабенко, О. В. Коваленко}

{\bf Об интерполяционных и экстремальных свойствах периодических идеальных сплайнов}
\end{center}
\end{Large}
\begin{abstract}
{Доказано существование и экстремальное свойство периодического идеального сплайна, интерполирующего заданную функцию в среднем.}
\end{abstract}

 Пусть $C^m$ ($m \in \ZZ_+$) обозначает пространство $m$ раз непрерывно дифференцируемых (непрерывных при $m=0$) $2\pi$-периодических  функций;  через $L_\infty$ будем обозначать пространство всех измеримых $2\pi$-периодичес\-ких функций с конечной нормой $\|f\| = \|f\|_\infty.$
Для $r\in\NN$ через $L^r_\infty$ будем обозначать пространство функций $f \in C$ таких, что функция $f^{(r-1)}$ абсолютно непрерывна, а  $f^{(r)}\in L_\infty$.

Для $r\in\NN$ функцию  $s(t)\in C^{r-1}$ назовем сплайном порядка $r$  с узлами в точках  
\begin{equation}\label{knots}
t_0<t_1<\ldots<t_n<t_{n+1}:=t_0 + 2\pi,
\end{equation} 
если на каждом из промежутков $[t_i,t_{i+1}]$, $i=0,1,\dots,n,$ $s(t)$ совпадает с сужением на этот промежуток некоторого алгебраического многочлена степени не выше $r$. 

Идеальным сплайном порядка $r$ назовем сплайн $s(t)$ порядка $r$ с узлами в точках~\eqref{knots}, y которого $s^{(r)}(t)=(-1)^i\epsilon$, $\epsilon\in\{1,-1\}$, $t\in (t_{i},t_{i+1})$, $i=0,\dots,n$.

Для $r,n\in\NN$ обозначим через $\Gamma^r_n$ множество $2\pi$--периодических идеальных сплайнов порядка $r$ с не более чем $n$ узлами.

Идеальные сплайны играют важную роль при точном решении многих экстремальных задач теории аппроксимации (см., например,~\cite{korn2} и~\cite{korn3}).

Основным результатом данной работы является следующая теорема.

\begin{theorem}\label{th1}
Пусть $r,m\in\NN$, $\varphi_k(x):\RR\to\RR$ --- интегрируемые четные функции с компактными носителями $[-\varepsilon_k,\varepsilon_k]$, положительные на $(-\varepsilon_k,\varepsilon_k)$, и такие, что 
$$
\int\limits_{-\varepsilon_k}^{\varepsilon_k}{\varphi_k(x)dx}=1,\,k=1,\dots,{2m+1}. 
$$

Пусть также числа $0<x_1<x_2<\ldots<x_{2m+1}<2\pi$ таковы, что носители $[x_k-\varepsilon_k,x_k+\varepsilon_k]$ функций $\varphi_k(x-x_k)$ попарно не пересекаются и содержатся в $[0,2\pi)$, $k=1,\dots,{2m+1}$. Тогда для любой функции $f\in L_{\infty}^r$ существуют число $\xi $ и сплайн $s(t)\in \xi\Gamma^r_{2m}$ такие,  что $$\int\limits_{x_k-\varepsilon_k}^{x_k+\varepsilon_k}{\varphi_k(x-x_k)} s(x)dx=\int\limits_{x_k-\varepsilon_k}^{x_k+\varepsilon_k}{\varphi_k(x-x_k)}f(x)dx ,\,{k=1,\dots,2m+1}.$$ При этом 
\begin{equation}\label{extremalProperty}
\|s^{(r)}\|_{\infty}=|\xi|\leq\|f^{(r)}\|_{\infty}.
\end{equation}

\end{theorem}

Из теоремы~\ref{th1}, в частности, следует возможность обыкновенной и кратной интерполяции периодическими идеальными сплайнами.

\begin{theorem}\label{th3} Пусть $r,m\in\NN$, $0\leq x_1<x_2<\ldots<x_{l}<2\pi$, $k_1,\dots,k_l \leq r-1$ -- целые неотрицательные числа такие, что $\sum_{j=1}^l{(k_j+1)}=2m+1$. Тогда для любой функции $f\in L_{\infty}^r$ существуют  число $\xi$ и сплайн $s(t)\in \xi\Gamma^r_{2m}$ такие,  что $$s^{(j)}(x_i)=f^{(j)}(x_i),\, j=0,\dots,k_i;\,i=1,\dots,l.$$ При этом выполняется экстремальное свойство~\eqref{extremalProperty}.

В частности, если заданы числа $0\leq x_1<x_2<\ldots<x_{2m+1}<2\pi$ и фукнция $f\in L^r_{\infty}$, то существуют число $\xi$ и сплайн  $s(t)\in \xi\Gamma^r_{2m}$ такие, что $$ s(x_i)=f(x_i),\, i=1,\dots,2m+1.$$ При этом выполняется экстремальное свойство~\eqref{extremalProperty}.
\end{theorem}

Теорема~\ref{th3} дает периодический аналог результатов, которые в непериодическом случае были получены С.~Карлиным~\cite{Karlin1, Karlin2}. Методы доказательства, используемые нами, значительно проще тех, которые применял С. Карлин, и могут быть использованы для решения других задач.

Приведем основные этапы доказательства теоремы~\ref{th1}. Теорема~\ref{th3} легко устанавливается с помощью теоремы~\ref{th1}.

Символом $*$ далее будем обозначать операцию свертки $2\pi$ - периодических функций. Пусть $f\in C$. Обозначим через $\nu(f)$ количество перемен знака $f$ на периоде. Для $\varepsilon>0$ и $x\in\RR$ положим
$$\A(x):=\frac1{2\pi}\sum_{j=-\infty}^{\infty}{\frac{{\rm e}^{ijx}}{\ch\varepsilon j}}.$$
Отметим некоторые свойства $A_\varepsilon$. Если $f\in C$, то $(\A*f)$ -- аналитическая функция на вещественной прямой (см.~\cite{Akhiezer37},\cite{Akhiezer38}; см. также~\cite{Nikolskiy46}, \S3); $(\A*f)$ равномерно сходится к $f$ при $\varepsilon\to 0$ на вещественной прямой; $\nu(\A*f)\le \nu(f)$ (см., напр.,~\cite{Mairhuber});
$\A$ -- четная функция; $\int_\RR{\A (x)dx}=1$. Всюду далее через $B_{r}(t)$, $r\in\NN$ мы будем обозначать ядро Бернулли порядка ${r}$ (см., например,~\cite{korn2}, \S 3.1).

Положим $C_k:=\int\limits_{x_k-\varepsilon_k}^{x_k+\varepsilon_k}{\varphi_k(x-x_k)f(x)}dx$, $k=1,2,\dots,2m+1$. Для $\varepsilon > 0$ рассмотрим экстремальную задачу 
\begin{equation}\label{extremalProblem}
\int\limits_0^{2\pi}{\left|c+\sum_{k=1}^{2m+1}{c_k(\A*B_r*\varphi_k)(x-x_k)}\right|dx}\to \min 
\end{equation}при условиях 
\begin{equation}\label{condition1}
\sum_{k=1}^{2m+1}{c_kC_k}=1,
\end{equation}
\begin{equation}\label{condition2}
\sum_{k=1}^{2m+1}{c_k}=0.
\end{equation}
Свертка с ядром $A_\varepsilon$ нам нужна для того, чтобы подынтегральная функция была аналитической, а значит минимизируемая функция была непрерывно дифференцируемой по параметрам $c$ и $c_k$, $k=1,\dots,2m+1$. Нетрудно доказать, что минимум в экстремальной задаче \eqref{extremalProblem} -- \eqref{condition2} существует. 

Введем функцию Лагранжа:
$$L:=\eta \int\limits_0^{2\pi}{\left|c+\sum_{k=1}^{2m+1}{c_k(\A*B_r*\varphi_k)(x-x_k)}\right|dx}+\lambda\left(\sum_{k=1}^{2m+1}{c_kC_k}-1\right)+\mu \sum_{k=1}^{2m+1}{c_k}.$$

Cогласно методу неопределенных множителей Лагранжа (см. например,~\cite{Metod_Lagr}, \S 2), существуют числа ${\eta_\varepsilon,\lambda_\varepsilon,\mu_\varepsilon\in\RR}$ такие, что $\eta_\varepsilon^2+\lambda_\varepsilon^2 + \mu_\varepsilon^2 \neq 0$ (мы можем, например, считать, что $\eta_\varepsilon^2+\lambda_\varepsilon^2 + \mu_\varepsilon^2 =1)$, и $c^\varepsilon,c_1^\varepsilon,c_2^\varepsilon,\dots,c_{2m+1}^\varepsilon\in\RR$, удовлетворяющие равенствам~\eqref{condition1} и~\eqref{condition2}  такие, что
\begin{equation*}
\eta_\varepsilon\int\limits_0^{2\pi}{{\rm sgn}\left[c^\varepsilon+\sum_{k=1}^{2m+1}{c_k^\varepsilon(\A*B_r*\varphi_k)(x-x_k)}\right]dx}=0
\end{equation*}
и для $j=1,2,\dots,2m+1$
\begin{equation*}
\eta_\varepsilon\int\limits_0^{2\pi}{{\rm sgn}\left[c^\varepsilon+\sum_{k=1}^{2m+1}{c_k^\varepsilon(\A*B_r*\varphi_k)(x-x_k)}\right](\A*B_r*\varphi_j)(x-x_j) dx}+\lambda_\varepsilon C_j+\mu_\varepsilon=0.
\end{equation*}

С помощью предельного перехода при $\varepsilon\to 0$ отсюда выводится, что существуют числа $\eta, \lambda, \mu$, $ \eta^2 + \lambda^2 + \mu^2 = 1$, и функция $g(x)$, почти всюду отличная от нуля и такая, что $\nu(g)\leq 2m$, для которых $\eta\int\limits_0^{2\pi}{{\rm sgn}g(x)}dx=0$,
и для $k=1,2,\dots,2m+1$
\begin{equation}\label{5''}
\eta\int\limits_0^{2\pi}{{\rm sgn}g(x)(B_r*\varphi_k)(x-x_k) dx}+\lambda C_k+\mu=0.
\end{equation}
Предположение о том, что $\eta = 0$ или $\lambda = 0$ приводит к противоречию. Это значит, что для функции $h(x):=({\rm sgn}g*B_r)(x)$ имеем  $h^{(r)}(x) ={\rm sgn} g(x)$, а значит $h(x)\in \Gamma^r_{2m}$.

Из \eqref{5''} получаем, что для $k=1,\dots,2m+1$
$$0=\eta\int\limits_0^{2\pi}{{\rm sgn}g(x)(B_r*\varphi_k)(x-x_k) dx}+\lambda C_k+\mu=$$ $$=(-1)^r\eta ({\rm sgn}g*B_r*\varphi_k)(x_k)+\lambda C_k+\mu=(-1)^r\eta(\varphi_k *h)(x_k)+\lambda C_k+\mu=$$ $$=(-1)^r\eta\int\limits_{x_k-\varepsilon_k}^{x_k+\varepsilon_k}{\varphi_k(x-x_k)}h(x)dx+\lambda C_k+\mu $$(последнее равенство верно в силу четности $\varphi_k(x)$).

Это значит, что идеальный сплайн  $s(x):=-\frac{(-1)^r\eta h(x)+\mu}{\lambda}$ является искомым интерполяционным сплайном.

Покажем, что $\|s^{(r)}\|_{\infty}\leq\|f^{(r)}\|_{\infty}$. Предположим противное. Пусть $\|s^{(r)}\|_{\infty}>\|f^{(r)}\|_{\infty}$. Положим $\delta(x):=s(x)-f(x)$. Тогда $\delta^{(r)}(x)=s^{(r)}(x)-f^{(r)}(x)$ может  менять знак только в узлах сплайна $s(x)$, а значит $\nu(\delta^{(r)})\le 2m$. Так как ${\int_{x_k-\varepsilon_k}^{x_k+\varepsilon_k}{\varphi_k(x-x_k)}\delta(x)dx=0}$, функция $\delta(x)$ не равна нулю тождественно ни на каком интревале  из $[0;2\pi)$ и $\varphi_k(x)$ -- неотрицательная функция $(k=1,\dots,2m+1)$, то на каждом интервале $(x_k-\varepsilon_k,x_k+\varepsilon_k)$ функция $\delta(x)$ меняет знак. Поскольку носители  функций $\varphi_k(x-x_k)$, $k=1,\dots,2m+1$ попарно не пересекаются, то $\nu(\delta)\geq 2m+1$, что в силу периодичности функции $\delta$ приводит к противоречию.
Теорема доказана.